\documentclass{amsart}

\newtheorem{theorem}{Theorem}[section]
\newtheorem{lemma}[theorem]{Lemma}

\newtheorem{corollary}[theorem]{Corollary}
\newtheorem{example}[theorem]{Example}

\newtheorem{remark}[theorem]{Remark}

\begin{document}

\title{Strictly convex space : Strong orthogonality and Conjugate diameters}
 
\author{Debmalya Sain, Kallol Paul and  Kanhaiya Jha}

\begin{abstract}
In a normed linear space X an element  x is said to be orthogonal to another element y  in the sense of Birkhoff-James, written as $ x \perp_{B}y, $  iff 
$ \| x \| \leq \| x + \lambda y \| $ for all  scalars $ \lambda.$ We prove that a normed linear space X is strictly convex iff  for any two elements x, y  of the unit sphere $ S_X$, $ x \perp_{B}y $  implies   $ \| x + \lambda y \| > 1~ \forall~ \lambda \neq 0. $ We apply this result to find a necessary and sufficient condition for a Hamel basis   to be a strongly orthonormal Hamel basis in the sense of Birkhoff-James in a finite dimensional real strictly convex space X. Applying the result we give an estimation for lower bounds of $ \| tx+(1-t)y\|, t \in [0,1] $ and $ \| y + \lambda x \|, ~\forall ~\lambda $ for all elements $ x,y \in S_X $ with $ x \perp_B y. $ We  find a necessary and sufficient condition for the existence of conjugate diameters through the points $ e_1,e_2 \in ~S_X $ in a real strictly convex  space of dimension 2. The concept of generalized conjuagte diameters  is then developed for a real strictly convex smooth space of finite dimension.

\end{abstract}
\maketitle


\section{Introduction.}
\noindent Suppose $ (X, \|.\|) $ is a normed linear space, real or complex. X is said to be strictly convex iff every element of the unit sphere $ S_X = \{ x \in X : \|x \| = 1 \} $ is an extreme point of the unit ball $ B_X = \{ x \in X : \|x \| \leq 1 \} .$ There are many equivalent characterizations of the strict convexity of a normed space, some of them given in [11,16] are\\
(i) If $ x,y \in S_X $ are linearly independent  then we have $ \| x + y \| < 2,$\\
(ii) Every non-zero continuous linear functional attains a maximum on at most one point of the unit sphere, \\
(iii) If $ \| x+y\| = \|x\| + \|y\|,~ x \neq 0 $ then $ y = c x $ for some $ c \geq 0. $ \\
\noindent A normed linear space X is said to be  uniformly convex  iff given $ \epsilon > 0 $ there exists $ \delta > 0 $ such that whenever
$ x, y \in S_X $ and $ \|x - y \| \geq \epsilon $ then $ \| \frac{x+y}{2} \| \leq 1 - \delta.$ \\ 
\noindent The number $ \delta(\epsilon) = \inf \{ 1 - \| \frac{x+y}{2} \|~  : ~ x,y \in S_X,~\|x-y\| \geq \epsilon \} $ is called the modulus of convexity of X. A space X is strictly convex iff $ \delta(2) = 1.$ \\
The concept of  strictly convex  and uniformly convex  spaces have been extremely useful in the studies of the geometry of Banach Spaces. One may go through  [2,5-8,10,11,14,16-18] for more information related to strictly convex spaces and [4,6,7,11] for uniformly convex spaces. \\

\noindent In a normed linear space X an element x is said to be orthogonal to another element y  in the sense of Birkhoff-James [3,10,11], written as $ x \perp_{B}y,  $  iff 
\[ \| x \| \leq \| x + \lambda y \| ~~ \mbox{for all  scalars} ~\lambda .\]
Recently in [13] we introduced the notions of strongly orthogonal set and strongly orthonormal Hamel basis in the sense of Birkhoff-James. For the sake of clarity we mention them briefly here.\\
\noindent \textbf{Strongly orthogonal in the sense of Birkhoff-James}: In a normed linear space X an element x is said to be strongly orthogonal to another element y  in the sense of Birkhoff-James, written as $ x \perp_{SB}y,  $  iff 
\[ \| x \| < \| x + \lambda y \| ~~ \mbox{for all  scalars} ~\lambda \neq 0 .\]
If $ x \perp_{SB}y  $ then $ x \perp_{B}y  $ but the converse is not true. In $ l^2(\infty) $ the element $ (1,0) $ is orthogonal to $(0,1)$ in the sense of Birkhoff-James but not strongly orthogonal.\\
\noindent \textbf{Strongly Orthogonal Set in the sense of Birkhoff-James}:  A finite set of elements $ \{x_1,x_2, \ldots ,x_n \} $ in a normed linear space X is said to be a strongly orthogonal set in the sense of Birkhoff-James   iff for each  $ i \in  \{ 1,2, \ldots , n \}, $
\[  \| x_{i} \| < \| x_i + \sum_{j=1,j\neq i}^{n} \lambda_j x_j \|,\]
whenever not all $ \lambda_j$'s  are 0. \\
In addition if $ \|x_i\| = 1, $ for each i, then the set is called  strongly orthonormal set in the sense of Birkhoff-James.\\
\noindent \textbf{ Strongly orthonormal Hamel basis in the sense of Birkhoff-James}: A finite set of elements $ \{x_1,x_2, \ldots , x_n \} $ in a normed linear space X is said to be a strongly orthonormal Hamel basis in the sense of Birkhoff-James iff the set is a Hamel basis of X and is a strongly orthonormal set in the sense of Birkhoff-James. As for example the set $ \{ (1,0,0, \ldots, 0), (0,1,0, \ldots, 0), \ldots (0,0, \ldots,0,1) \} $ is a strongly orthonormal Hamel basis in the sense of Birkhoff-James in $ l^n(p) $ for $(1 \leq p < \infty) $  but not in $ l^n(\infty).$\\

\noindent We show that a normed linear space X is strictly convex iff  for $ x, y \in  S_X $ we have $ x \perp_{B}y $ implies  $ x \perp_{SB}y. $  Using this result we find a necessary and sufficient condition for a Hamel basis to be a strongly orthonormal Hamel basis in the sense of Birkhoff-James. We also show that in a normed linear space X if $ x, y \in  S_X $ and  $ x \perp_{B}y $ then $ \| y + \lambda x \| \geq \frac{1}{2}, $ for all $ \lambda.$ If the space is strictly convex then we  show that $ \| y + \lambda x \| > \frac{1}{2}$ for all $ \lambda.$ That the condition is only necessary but not sufficient is illustrated with an example. The result does not get better in uniformly convex spaces too. Following Theorem 1  of James[9] we conclude that a real normed linear space X (of dimension $ > 2 $) is  an inner product space iff for $ x, y \in  S_X $ we have $ x \perp_{B}y $ implies   $ y \perp_{SB}x.  $  We also find a necessary and sufficient condition for two diameters of $S_X$ to be conjugate diameters  in a  real strictly convex  space of dimension 2.  Thus given any two diameters in a strictly convex space we can say whether they are conjugate or not.  Finally we show that if  X is a real strictly convex smooth space of dimension 2 then $S_X$ is a Radon curve iff at each point of $ S_X$ there exists a strongly orthonormal Hamel basis  in the sense of Birkhoff-James containing that point. The concept of generalized conjugate diameters  is then developed for a real strictly convex  smooth space of finite dimension higher than 2.

 \section{Main Results}
 \noindent We first obtain the characterization of strictly convex spaces. 
\begin{theorem}
Suppose X is a real normed linear space. If for $ x, y \in X - \{0\}  $,  $ x \perp_{B} y $ implies   $ x \perp_{SB}y  $ then X is strictly convex.
\end{theorem}
\noindent \textbf{Proof.} If possible let the unit sphere $ S_X $ contains a straight line segment i.e., there exists $ \|u \| = \| v \| = 1 $ with $ \| tu + (1-t)v \| = 1 ~~ \forall t \in [0,1].$ \\
Let $ x = \frac{1}{2} u + \frac{1}{2} v,~ y = v - u .$ Consider $ x + \lambda y.$\\
If $ - \frac{1}{2} \leq \lambda \leq \frac{1}{2} $ then $ \| x + \lambda y \| = \| t u + (1-t) v \| = 1 ~~where ~~ t = \frac{1}{2} - \lambda. $\\
If $ \lambda < - \frac{1}{2} $ then $ \frac{1}{2} - \lambda > 1 $ and so we can write $ \frac{1}{2} - \lambda = t \alpha $ for some $ t \in (0,1)$ and $ \alpha > 1. $ In this case 
\begin{eqnarray*}
 \| x + \lambda y \| & = & \| t \alpha~ u + (1-t\alpha) ~ v \| \\
                     & = & \| t \alpha~ u + (\alpha -t\alpha) ~ v  + ( 1 - \alpha ) v \|\\
                     & \geq & \mid \| t \alpha~ u + (\alpha -t\alpha) ~ v \| - \| ( \alpha - 1) v \| \mid \\
                     & = & \mid ~~\mid \alpha \mid - \mid \alpha - 1 \mid ~~\mid \\
                     & = & 1
                     \end{eqnarray*}                     
\noindent If $ \lambda > \frac{1}{2} $ then $ \frac{1}{2} +\lambda > 1 $ and so we can write $ \frac{1}{2} + \lambda = t \alpha $ for some $ t \in (0,1)$ and $ \alpha > 1. $ In this case 
\begin{eqnarray*}
 \| x + \lambda y \| & = & \| (1- t \alpha)~ u + t\alpha ~ v \| \\
                     & = & \| ( \alpha - t \alpha)~ u + t\alpha ~ v  + ( 1 - \alpha ) u \|\\
                     & \geq & \mid \| ( \alpha -t \alpha)~ u + t\alpha ~ v \| - \| ( \alpha - 1) u \| \mid \\
                     & = & \mid ~~\mid \alpha \mid - \mid \alpha - 1 \mid ~~\mid \\
                     & = & 1
                     \end{eqnarray*}             
\noindent    Thus $ \| x + \lambda y \| \geq \|x\|~ \forall \lambda $ but $ \| x + \lambda_0 y \| = \|x \|,~ for~ ~\lambda_0 \in [-1/2, 1/2].$                 
This is a contradiction to our hypothesis and so X is strictly convex.\\
Alternatively the result can be proved using the convexity of the unit ball $ B_X$ and noting that  for each $ \lambda, ~ x + \lambda y = \alpha u + \beta v$ where  $ \alpha + \beta = 1. $ \\

\noindent Conversely we prove the following
\begin{theorem}
Suppose X is a strictly convex  space and  $ x, y \in X - \{0\} $ with $ x \perp_{B} y, $  then  $ x \perp_{SB} y. $ 
\end{theorem}
\noindent \textbf{Proof.} Without loss of generality we assume that there exists $ x, y \in X,~ \| x\| = 1$ such that $ \| x + \lambda y \| \geq 1~ \forall~ \lambda $ but
$ \| x + \lambda_0 y \| =  1 $ for some $ \lambda_0 \neq 0 .$\\
Let $ 0 < t < 1,$ then 
\begin{eqnarray*}
1 = t \|x \| + (1-t) \| x + \lambda_0 y \| & \geq & \| tx + (1-t) ( x + \lambda_0 y) \| \\
\Rightarrow 1 \geq \| tx + (1-t) ( x + \lambda_0 y) \| & = & \| x + (1-t) \lambda_0 y \| \geq 1
\end{eqnarray*}
Thus $ \| tx + (1-t) ( x + \lambda_0 y) \| = 1$ --
which contradicts the fact that X is strictly convex.

\begin{corollary}
If a normed linear space X is strictly convex and $ x_1,x_2,\ldots , x_n \in S_X ~ (n \geq 2) $ are linearly independent then we have $ \| x_1 + \lambda_2 x_2 + \ldots + \lambda_n x_n \| \geq 1 $ for all $ \lambda_2, \lambda_3, \ldots ,\lambda_n $ implies $ \| x_1 + \lambda_2 x_2 + \ldots + \lambda_n x_n \| > 1 $ for all $ \lambda_2, \lambda_3, \ldots ,\lambda_n $ with $ (\lambda_2, \lambda_3, \ldots ,\lambda_n) \neq (0,0, \ldots,0). $
\end{corollary}
\noindent \textbf{Proof.} Follows easily from the last theorem.\\

\noindent Thus we have the following characterization of a real strictly convex space which follows easily from the last two theorems.
\begin{theorem}
A real normed linear space X is strictly convex iff 
\[   x,y \in S_X ~ \mbox{and}~ x \perp_{B}y   \Rightarrow   x \perp_{SB} y.\]
\end{theorem}

\noindent James[9] proved that a real normed linear space X of   dimension $> ~2$  is an inner product space iff Birkhoff-James orthogonality is symmetric.
The following characterization of a real inner product space now follows easily from  Theorem 2.4 mentioned above and Theorem 1 of James[9].

\begin{theorem}
A real normed linear space X of   dimension $> ~2$  is an inner product space iff   
\[ x,y \in S_X ~ \mbox{and}~ x \perp_{B}y  \Rightarrow   x \perp_{SB} y  ~\mbox{and}~y \perp_{SB} x.  \]
\end{theorem}

\noindent We next try to calculate a lower bound for $ \| y + \lambda x \| ~~\forall~ \lambda  $ if $ x,y \in S_X ~ \mbox{and}~ x \perp_{B}y $ in a  normed linear space. We first prove the following two theorems: 

\begin{theorem}
Suppose X is a normed linear space and $ x,y \in S_{X} $ with $ \| x + \lambda y \| \geq 1~ \forall~ \lambda.$ Then there exists a $ \delta \geq \frac{1}{3} $ such that $ \| tx + (1-t)y \| \geq \delta~ \forall t \in [0,1].$  In addition $ \delta > \frac{1}{3} $ if X is  strictly convex.
\end{theorem}
\noindent \textbf{Proof.} Consider $ f : [0,1] \longrightarrow [0,1] $ defined as $ f(t) = \| tx + (1-t)y \| .$ Then f is a continuous function  on a compact set and so attains its minimum value. As $ f(0) = 1 $  by the property of continuity we can find $ \delta_0 \in (0,1) $ such that $ f(t) > \frac{1}{2}~ \forall ~t \in (0, \delta_0).$ \\
Again for all t with  $  t   \geq \delta_0 $ we get $ f(t) = \| tx + (1-t)y \| \geq \delta_0. $ Thus we get 
\[ \| tx + (1-t)y \| \geq \min \{\delta_0, \frac{1}{2} \} = \delta ~ \forall t \in [0,1].\] 

\noindent Clearly for all $ t \in [0,1],~~ \| tx + (1-t)y \| \geq \mid t \mid. $ \\
Also for all $ t \in [0,1],~~ \| tx + (1-t)y \| \geq \mid \mid t \mid - \mid (1-t) \mid \mid = \mid 2t-1 \mid.$
Thus $ \| tx + (1-t)y \| \geq \frac{1}{3}~ \forall t \in [0,1].$ \\
If X is strictly convex then by theorem $ 2.2,$ $ \| tx + (1-t)y \| > t $ for all $ t \in [0,1)$, also for all $ t \in [0,1],~~ \| tx + (1-t)y \| \geq \mid 2t-1 \mid $ and so 
$  \delta > \frac{1}{3}. $ 

\begin{remark}
Geometrically this shows that if an element x of the unit sphere is orthogonal to an element y of the unit sphere then the line segment joining x and y lies outside the open ball of radius $ \frac{1}{3}.$
\end{remark}

\begin{remark}
 If X is finite dimensional then there exists $ \delta > 0 $ such that $ \| tx + (1-t)y \| \geq \delta~ \forall t \in [0,1]$ and for all $ x,y \in S_{X} $ with $ x \perp_{B}y$ ( consider the continuous function $ f : M \times [0,1] \longrightarrow R $  defined on the compact set $ M \times [0,1]  $ as $ f(x,y,t) = \| tx + (1-t)y \| $  where $ M = \{ (x,y) \in S_X \times S_X : x \perp_{B} y \}$ ). 
\end{remark}
\begin{theorem}
Suppose $ (X,\| .\|) $ is a normed linear space and  $ x,y \in S_{X}$  such  that $  \|x+ \lambda y\| \geq 1~ \forall \lambda .$ Then there exists  $ \delta \geq \frac{1}{2}  $ such  that $  \|y+ \lambda x\| \geq \delta~ \forall \lambda.$ In addition $ \delta > \frac{1}{2} $  if X is a  strictly convex normed linear space.
\end{theorem} 
\noindent \textbf{Proof.} Clearly for all $ \lambda , ~~\| y + \lambda x \| \geq \mid \lambda \mid.$ \\
Also for all $ \lambda, ~~\| y + \lambda x \| \geq \mid  1 - \mid \lambda \mid \mid.$ \\
Thus for all $ \lambda , ~~\| y + \lambda x \| \geq \frac{1}{2}.$  \\
\noindent If X is strictly convex then using theorem $ 2.1 $   , $ \| y + \lambda x \| > \mid \lambda \mid $ for all $ \lambda $  and  $~~\| y + \lambda x \| \geq \mid  1 - \mid \lambda \mid \mid$  so that $ \delta > \frac{1}{2} $.

\begin{remark}
Geometrically this shows that if an element x of the unit sphere is orthogonal to an element y of the unit sphere then the line through  y parallel to x  lies outside the open ball of radius $ \frac{1}{2}.$
\end{remark}

\begin{remark}
If X is finite dimensional then one can find a $ \delta > 0 $ such that $ \| y + \lambda x \| \geq \delta~ \forall \lambda $ and for all $ x,y \in S_{X} $ with $ x \perp_{B}y$ ( consider the continuous function $ f : M \times S_K \longrightarrow R $  defined on the compact set $ M \times S_K $ as $ f(x,y,\lambda) = \| y + \lambda x \| $  where $ M = \{ (x,y) \in S_X \times S_X : x \perp_{B} y \} $ ).
\end{remark}

\noindent We next give an example to show that the bounds can not be improved on in an arbitrary normed linear space.
\begin{example}
 Consider $(R^2,\| . \|_\infty)$  with $ x =(1,1) $ and $ y = (-1,0) $, then  $ \| x + \lambda y \| \geq \|x\|~ \forall~ \lambda $ but
 $ \|tx+(1-t)y\| = \frac{1}{3} $ for  $ t = \frac{1}{3}. $  Also   $ \|y+ \lambda x \| = \frac{1}{2}  $ for $ \lambda = \frac{1}{2}.$ This shows that in an arbitrary normed linear space one can not improve on the bounds.
 \end{example} 
\noindent That the condition of $ \delta > \frac{1}{2} $ is necessary but not sufficient for strict convexity  follows from the example given below.
\begin{example}
Consider the normed linear space X with 
$ S_X = \{ (x,y) : -1 \leq x \leq 1,~ y = \pm 1 \} ~ \cup ~ \{ (x,y): 1 \leq x \leq \sqrt{2}, y = \pm \sqrt{2 - x^2} \} ~ \cup ~ \{ (x,y): - \sqrt{2} \leq x \leq -1, y = \pm  \sqrt{2 - x^2} \}.$  Then it is easy to see that X is not strictly convex and with a little bit of calculation we can show that if $ x,y \in S_{X}$  such  that $  \|x+ \lambda y\| \geq \|x\|~ \forall \lambda $ then there exists  $ \delta >\frac{1}{2} $ such  that $  \|y+ \lambda x\| > \delta~ \forall \lambda.$
\end{example}
\noindent We can sum up our results in  the following two theorems which give a better idea of how the estimation of the bounds discussed above behave in different types of normed linear spaces.

\begin{theorem}
Suppose X is a normed linear space. 
\begin{itemize}
	\item If $ x, y \in S_X $ and $ x \perp_{B} y $ then $ \| tx + (1-t)y \| \geq \frac{1}{3} ~~\forall t \in [0,1].$
	\item If X is strictly convex then  $ x, y \in S_X $ and $ x \perp_{B} y $ implies $ \| tx + (1-t)y \| > \frac{1}{3} ~~\forall t \in [0,1] .$
	\item If X is an inner product space then $ x, y \in S_X $ and $ x \perp_{B} y $ implies $ \| tx + (1-t)y \| \geq \frac{1}{\sqrt{2}}~ ~\forall t \in [0,1]. $
\end{itemize}
\end{theorem}

\begin{theorem}
Suppose X is a normed linear space. 
\begin{itemize}
	\item If $ x, y \in S_X $ and $ x \perp_{B} y $ then $ \| y + \lambda x \| \geq \frac{1}{2} ~\forall ~\lambda.$
	\item If X is strictly convex then  $ x, y \in S_X $ and $ x \perp_{B} y $ implies $ \| y + \lambda x \| > \frac{1}{2} ~\forall ~\lambda .$
	\item If X is an inner product space then $ x, y \in S_X $ and $ x \perp_{B} y $ implies $ \| y + \lambda x \| > 1~~ (\lambda \neq 0). $
\end{itemize}
\end{theorem}
\begin{remark}
 Considering the uniformly convex spaces $ (R^2,l_p)~(p > 1) $ it may be remarked that  a better estimation of the bounds  discussed above may not be possible in an arbitrary uniformly convex space.
\end{remark}

\section{ Strongly orthonormal Hamel basis in the sense of Birkhoff-James}
\noindent In this section we find a necessary and sufficient condition for a Hamel basis to be a strongly  orthonormal Hamel basis in the sense of Birkhoff-James in a finite dimensional real strictly convex space. We begin with the following lemma.
\begin{lemma}
Suppose X is a finite dimensional real normed linear space and  $ \{ e_1,e_2, \ldots , e_m\} $ be a Hamel basis of X with $ \|e_i\| = 1, \forall i =1,2, \ldots , m.$ \\
For $ i =1,2, \ldots ,m $ let
\[ S_i = \{ \alpha_{iz} ~: ~ z = \alpha_{1z}e_1 + \alpha_{2z}e_2 + \ldots + \alpha_{mz}e_m \in S_X~\}. \] 
Then $S_i $ attains its maximum value for each $ i =1,2, \ldots , m.$ 
\end{lemma}
\noindent \textbf{Proof.} We first show that each $ S_i $ is closed and bounded.\\
Let $ \{ \alpha_{i n} \}_{n=1}^{\infty} $ be a sequence in $ S_i $ which converges to some real number $ \beta $  where $ z_n = \alpha_{1n}e_1 + \alpha_{2n}e_2 + \ldots + \alpha_{mn}e_m $ and $ \|z_n \| = 1 $ for all n.\\
Since $ S_X $ is compact so $ \{ z_n\} $ has a convergent subsequence $ z_{n_k} $ converging to some $ z_0 $ in $ S_X.$ 
Let $ z_0 = \gamma_1 e_1 + \gamma_2 e_2 + \ldots + \gamma_m e_m $. Then 
\begin{eqnarray*}
\| z_{n_k} - z_0 \| & = &  \mid (\alpha_{1n_k} - \gamma_1) e_1 + (\alpha_{2n_k} - \gamma_2) e_2 + \ldots + (\alpha_{mn_k} - \gamma_m) e_m \mid \\
 & > & c [~ \mid (\alpha_{1n_k} - \gamma_1)  \mid  + \mid (\alpha_{2n_k} - \gamma_2) \mid + \ldots + \mid (\alpha_{mn_k} - \gamma_m) \mid~ ] 
 \end{eqnarray*}
for some $ c > 0 $ as any two norms on a finite dimensional space are equivalent.\\
As $ \{z_{n_k}\} $ converges to $ z_0 $ so $ \{ \alpha_{in_k} \} $ converges to $ \gamma_i$ as $ k \longrightarrow \infty. $\\
This shows that $ \beta = \gamma_i.$ Thus $ S_i $ is closed for each i. \\
Similarly using the fact that  any two norms on a finite dimensional space are equivalent it is easy to show that each $ S_i $ is bounded. So  each $ S_i $ is compact and hence attains its maximum value. This completes the proof of the lemma.\\

\noindent In the next lemma we show that if X is strictly convex then there exists a  unique element of $ S_X $ for which  the coefficient of $ e_i $ attains the maximum value of $ S_i. $ 
\begin{lemma}
Suppose X is a finite dimensional  real strictly convex  space. Let $ \alpha_{i_0z_{0}} = \max S_{i_0} $ which is attained for $ z_0 = \alpha_{1z_0}e_1 + \alpha_{2z_0}e_2 + \ldots + \alpha_{i_0z_{0}}e_{i_0} +\ldots + \alpha_{mz_0}e_m \in S_X.$ If $ y = \beta_1e_1 + \beta_2e_2 + \ldots + \beta_m e_m \in S_X $ with $ \beta_{i_0} = \alpha_{i_0z_{0}} $ then $ y = z_0.$
\end{lemma}
\noindent \textbf{Proof.} For any $ t\in [0,1] $ we have $ \|ty+(1-t)z_{0}\| \leq 1. $ \\
Also for any $ t\in [0,1] $, 
\begin{eqnarray*}
 ty+(1-t)z_{0} & = & t(\beta_1e_1 + \beta_2e_2 + \ldots + \beta_m e_m)+(1-t)~\\
  & & (~\alpha_{1z_0}e_1 + \alpha_{2z_0}e_2 + \ldots + \alpha_{i_0z_{0}}e_{i_0} + \ldots + \alpha_{mz_0}e_m~) \\
 & = & \gamma_1 e_1 + \gamma_2 e_2 + \ldots +  \gamma_{i_0} e_{i_{0}} + \ldots + \gamma_m e_m
 \end{eqnarray*}
 where $ \gamma_1 = t\beta_1+(1-t)\alpha_{1z_0}, \gamma_2 = t \beta_2+(1-t)\alpha_{2z_0}, \ldots, \gamma_{i_0} = \alpha_{i_0z_{0}}, \ldots ,\gamma_m = t\beta_m+(1-t)\alpha_{mz_0}. $ \\
 If $ \|ty+(1-t)z_{0}\| < 1 $ then $ \frac{ ty+(1-t)z_{0}} { \|ty+(1-t)z_{0}\| } \in S_X $ and 
 \[ \frac{ ty+(1-t)z_{0}} { \|ty+(1-t)z_{0}\| } = \frac{1}{\|ty+(1-t)z_{0}\|} (~\gamma_1 e_1 + \gamma_2 e_2 + \ldots +  \gamma_{i_0} e_{i_{0}} + \ldots + \gamma_m e_m~) .\]
 As $ \frac{\gamma_{i_0}}{\|ty+(1-t)z_{0}\|} > \gamma_{i_0} = \alpha_{i_0z_{0}} = \max S_{i_0} $ so we get a contradiction.\\ 
 Hence $ \|ty+(1-t)z_{0}\| = 1 $ for all $ t \in  [0,1]. $ \\
 Since $(X,\|.\|) $  is a strictly convex space so we must have $ y = z_0. $ \\

\noindent We finally prove the theorem which gives the necessary and sufficient condition for a Hamel basis to be  a strongly orthonormal Hamel basis in the sense of Birkhoff-James.
\begin{theorem}
Suppose X is a finite dimensional real strictly convex  space and $ \{ e_1,e_2, \ldots , e_m\} $ be a Hamel basis of X with $ \|e_i\| = 1, \forall i =1,2, \ldots , m.$  Then a necessary and sufficient condition for  $ \{ e_1,e_2, \ldots , e_m\} $ to be a strongly  orthonormal Hamel basis in the sense of Birkhoff-James is 
$  \max S_i = 1 ~\forall~ i = 1,2, \ldots ,m. $
\end{theorem}
\noindent \textbf{Proof. } Let $ \max S_i = 1 $ for all $ i = 1,2, \ldots , m. $\\
Now 
\[ \| e_i + \sum_{j=1,j\neq i}^{m} \lambda_j e_j \| \geq 1~ \forall~ \lambda_j~ \mbox{with} ~ 1 \leq j \leq m,~ j \neq i. \]
for otherwise the coefficient of $ e_i $ in the unit element $ \frac{e_i + \sum_{j=1,j\neq i}^{m} \lambda_j e_j}{\| e_i + \sum_{j=1,j\neq i}^{m} \lambda_j e_j \|} $ will be greater than 1 which is not possible as $  \max S_i = 1. $
Since the normed space $ (X,\|.\|) $ is strictly convex so by applying Corollary $ 2.3 ,$ we obtain that 
 \[ \| e_i + \sum_{j=1,j\neq i}^{m} \lambda_j e_j \| > 1~ \forall~ \lambda_j~ \mbox{with} ~ 1 \leq j \leq m,~ j \neq i,~ \mbox{whenever not all} ~\lambda_{j}'s~ \mbox{are}~ 0.\] 
\noindent This proves that $ \{ e_1,e_2, \ldots , e_m\} $ is a strongly orthonormal Hamel basis in the sense of Birkhoff-James  of  $ (X,\|.\|).$ \\
\noindent Conversely let $ \{ e_1,e_2, \ldots , e_m\} $ be a strongly orthonormal Hamel basis in the sense of Birkhoff-James.\\
For ~$  z=\alpha_{1z}e_{1}+ \alpha_{2z}e_{2}+ \ldots + \alpha_{mz}e_{m} \in S_{X} $, we have \\
$ 1 = \|z\| = \|\alpha_{1z}e_{1} + \alpha_{2z}e_{2} + \ldots + \alpha_{mz}e_{m} \| \geq |\alpha_{iz}| ~\forall~ i = 1,2,\ldots,m.$ \\
Thus $  \max ~S_{i} \leq 1.$ On the other hand, $ \max ~S_{i} \geq 1 $ since $ \|e_{i}\|=1.$
This proves that $ \max ~S_{i} = 1 ~  \forall~ i = 1,2,\ldots,m. $ \\

\begin{remark}
In $ (R^2, l_p),(p>1) $ the set $ \{ (1,0),(0,1) \} $ is a strongly orthonormal Hamel basis in the sense of Birkhoff-James since the maximum possible value of the $x$-coordinate of any point on the unit sphere of  $ (R^2, l_p) $ is attained at the point (1,0),  the same holds for the point (0,1).
\end{remark}
\section{Conjugate diameters and strongly orthonormal Hamel basis in the sense of Birkhoff-James} 
\noindent If the dimension of the real strictly convex  space X is 2 then the concept of  strongly orthonormal Hamel basis in the sense of Birkhoff-James is  connected with the concept of conjugate diameters. We first give the definition of conjugate diameters in a conic.  Two diameters of a conic are said to be conjugate iff each of them bisects the chords of the conic that are parallel to the other. If X is a real normed linear space with $ dim~X = 2 $ then  two diameters of the unit sphere $ S_X$ passing through $x$ and $y$ are said to be conjugate  iff $ x\bot_{B}y$ and $ y \bot_{B}x.$ It follows from Auerbach [1] that $ S_X$ has at least one pair of conjugate diameters. We now state the following theorem the proof of which is clear from definition of conjugate diameters and Theorem 2.3.
\begin{theorem}
Suppose X is a real strictly convex   space of dimension 2. Then $S_X$ has a pair of conjugate diameters through $ e_1,e_2 $ iff  $ \{ e_1,e_2 \}$ is a strongly orthonormal Hamel basis  in the sense of Birkhoff-James.
\end{theorem}
\noindent We now give an example of a real normed linear space of dimension 2 which has no strongly orthonormal Hamel basis in the sense of Birkhoff-James but has a pair of conjugate diameters to show that strict convexity in the last theorem is necessary. 
\begin{example}
Consider $(R^2,\|.\|)$ where the unit sphere is given by $ S_{R^2} = \pm A \cup \pm B \cup \pm C$ where \\
$ A  =  \{ (x,y):~ x^4 + y^4 = 1, ~ 0 \leq x \leq \frac{1}{3},~ y \geq 0 \} $ \\
$  B  =  \{ (x,y):~ (-x)^3 + y^3 = 1, ~ - \frac{1}{4} \leq x \leq 0,~ y \geq 0 \}  $ \\
 $ C  =  \{ t(\frac{1}{3},y_1) + (1-t)(\frac{1}{4},y_2):~ 0 \leq t \leq 1, (\frac{1}{3})^4 + (y_1)^4 = 1, (\frac{1}{4})^3 + (-y_2)^3 = 1 \}.$ \\ 
 \noindent It is easy to check that $(R^2,\|.\|)$ has no strongly orthonormal Hamel basis in the sense of Birkhoff-James but has a pair of conjugate diameters.
\end{example}
\noindent We next explore the relation between the concept of strongly orthonormal Hamel basis  in the sense of Birkhoff-James and the notion of Radon curve [12,15], the definition of which is given below.\\
\noindent \textbf{Radon Curve.} A curve is a Radon curve iff it is the unit sphere of a two dimensional real normed linear space (also known as Minkowski plane) in which the notion of orthogonality (also called normality) in the sense of Birkhoff-James is symmetric. We now give the following theorem.
\begin{theorem}
Suppose X is a real strictly convex  smooth space of dimension 2. Then $S_X$ is a Radon curve iff at each point of $ S_X$ there exists a strongly  orthonormal Hamel basis  in the sense of Birkhoff-James containing that point.
\end{theorem}
\noindent{\textbf{Proof.}} The proof follows from the fact that in a smooth space the supporting line to $S_X$ at any point on $S_X$ is unique.\\

\noindent So using the concept of strongly orthonormal Hamel basis  in the sense of Birkhoff-James  one can generalize the notion of conjugate diameters and Radon curve in a real strictly convex  smooth space even if the dimension of the space is greater than 2. \\

\noindent\textbf{Generalized conjugate diameters.} Suppose X is a real normed linear space of dimension n. Then a set of n diameters of $ S_X$ are said to be generalized conjugate diameters iff any two of them  are conjugate diameters.  
As for example in $(R^n,l_p)~(1<p< \infty)$ diameters passing through $ (1,0,\ldots,0), (0,1,0, \ldots,0),$ $ \ldots \ldots , (0,0,\ldots,1) $ are generalized conjugate diameters since \\$ \{ (1,0,\ldots,0),  (0,1,0, \ldots,0), \ldots \ldots ,(0,0,\ldots,1) \} $ is a strongly orthonormal Hamel basis  in the sense of Birkhoff-James in $(R^n,l_p).$\\

\section{Acknowledgements}
  The first author would like to thank Jadavpur University and DST, Govt. of India  for the partial financial support provided through DST-PURSE project and the second author would like to thank UGC, Govt. of India for the financial support. \\
  We would like to thank Professor T. K. Mukherjee for his invaluable suggestion while preparing this paper.




\bigskip
\noindent
\parbox[t]{.48\textwidth}{
Debmalya Sain and Kallol Paul\\
Department of Mathematics\\
Jadavpur University \\
Kolkata 700032\\
INDIA\\
kalloldada@gmail.com\\
 saindebmalya@gmail.com } \hfill
\parbox[t]{.48\textwidth}{
Kanhaiya Jha\\
Department of Mathematical Sciences\\
 School of Science, Kathmandu University,\\ 
POBox Number 6250, Kathmandu, NEPAL
}

\end{document}